\documentclass{gtpart}   
\usepackage{pinlabel}
\usepackage[all]{xy}
\usepackage{amsmath,amsthm}
\usepackage{amssymb}
\usepackage[all]{xypic}

\title{Realizing Families of\\
Landweber Exact Homology Theories}

\author{Paul G. Goerss}
\givenname{Paul}
\surname{Goerss}
\address{Department of Mathematics\\Northwestern University\\\newline
         Evanston, IL 60208\\USA}
\email{pgoerss@math.northwestern.edu}
\urladdr{http://www.math.northwestern.edu/~pgoerss}


\swapnumbers
\newtheorem{thm}{Theorem}[section]

\newtheorem{prop}[thm]{Proposition}

\newtheorem{defn}[thm]{Definition}
\theoremstyle{definition}
\newtheorem{rem}[thm]{Remark}

\newtheorem{exam}[thm]{Example}

\newtheorem{warn}[thm]{Warning}
\newtheorem{realize}[thm]{The Realization Problem}
\numberwithin{equation}{section}
\numberwithin{figure}{section}

\def\Aut{\operatorname{Aut}}

\def\ZZ{{{\mathbb{Z}}}}
\def\PP{{{\mathbb{P}}}}
\def\QQ{{{\mathbb{Q}}}}

\def\MR{{{\mathcal{TM}}}}

\def\cE{{\cal E}}
\def\cF{{\cal F}}

\def\cH{{\cal H}}

\def\cM{{\cal M}}
\def\cN{{\cal N}}
\def\cO{{\cal O}}

\def\cU{{\cal U}}

\def\la{{\lambda}}

\def\colim{\operatorname{colim}}
\def\holim{\operatorname{holim}}

\def\Ext{{{\mathrm{Ext}}}}
\def\Ga{\Gamma}

\def\Hom{\mathrm{Hom}}

\def\map{\operatorname{map}}

\def\longr{{\ \longrightarrow\ }}
\def\defeq{\overset{\mathrm{def}}=}
\def\mm{{{\mathfrak{m}}}}

\def\Spec{{{\mathrm{Spec}}}}
\def\GG{{{\mathbb{G}}}}
\def\FF{{{\mathbb{F}}}}
\def\bFF{{{\bar{\FF}_p}}}

\def\MR{{{\mathcal{TM}}}}

\def\fg{{\mathbf{fg}}}
\def\fgl{{\mathbf{fgl}}}

\def\cN{{{\mathcal{N}}}}
\def\Spf{{{\mathrm{Spf}}}}

\def\Def{{{\mathbf{Def}}}}

\date{}

\def\CP{{{\mathbb{C}\mathrm{P}}}}
\def\for{{{\mathrm{for}}}}
\def\Ell{{{e\ell\ell}}}
\def\mweier{{{\cM_{\mathrm{Weier}}}}}
\def\mell{{{\cM_\Ell}}}

\def\\tmf{{{\mathbf{\tmf}}}}

\def\Top{{{\mathbf{top}}}}

\def\tmf{{{\mathbf{tmf}}}}

\def\barell{{{\bar{\cM}_\Ell}}}

\begin{document}

\begin{abstract} I discuss the problem of realizing families
of complex orientable homology theories as families of $E_\infty$-ring
spectra, including a recent result of Jacob Lurie emphasizing the role of
$p$-divisible groups.\footnote{The author was partially supported by the National
Science Foundation.}
\end{abstract}
\maketitle

A few years ago, I wrote a paper \cite{Newton} discussing a realization
problem for families of Landweber exact spectra. Since Jacob Lurie \cite{Lurie}
now has a major positive result in this direction, it seems worthwhile
to revisit these ideas. 

In brief, the realization problem can be stated as follows. Suppose we
are given a flat morphism
$$
g:\Spec(R) \longr \cM_\fg
$$
from an affine scheme to the moduli stack of smooth 1-dimensional formal groups. 
Then we get a two-periodic
homology theory $E(R,G)$ with $E(R,G)_0 \cong R$ and associated formal
group
$$
G = \Spf(E^0\CP^\infty)
$$
isomorphic to the formal group classified by $g$.
The higher homotopy groups of $E(R,G)$ are zero in odd degrees and
$$
E(R,G)_{2n} \cong \omega_G^{\otimes n}
$$
where $\omega_G$ is the module of invariant differentials for $G$. The
module $\omega_G$
is locally free of rank $1$ over $R$, and free of rank 1 if $G$ has a coordinate.
In this case $E(R,G)_\ast = R[u^{\pm 1}]$ where $u \in E(R,G)_2$ is
a generator.
The fact that $g$ was flat implies $E(R,G)$ is Landweber exact,
even if $G$ doesn't have a coordinate.

Now suppose we are given a flat morphism of stacks
$$
g:X \longr \cM_\fg.
$$
Then for each affine open $\Spec(R) \to X$, we get a spectrum $E(R,G)$
and, because there are no phantom maps between these spectra 
by \cite{HS}, we get a presheaf $\cO^\top_X$ on $X$ in the homotopy category
with
$$
\cO^\top_X(\Spec(R) \to X) = E(R,G).
$$

A na\"\i ve version of the realization problem is this: can this presheaf
$\cE_X$ be lifted to a presheaf (or sheaf) $\cO^\Top_X$ of $E_\infty$-ring
spectra? If so, how unique is this lift?

I call this na\"\i ve because, at the very least, we want some hypotheses
on $X$ and the morphism $g:X \longr \cM_\fg$. For example, we might want
to specify that $X$ be an algebraic stack, or a Deligne-Mumford stack,
and we might want to specify that the morphism $g$ be representable.
Even having done this, I don't suppose anyone expects a positive
answer in this generality -- there are simply too many flat maps
to $\cM_\fg$. Indeed, Lurie's result \ref{mainthm}
below requires that $g$ factor as
$$
\xymatrix{
X \rto^-f &\cM_p(n) \rto & \cM_\fg
}
$$
where $\cM_p(n)$ is a moduli stack of $p$-divisible groups and
$f$ is appropriately \'etale. As a consequence, we don't just have
a family of formal groups over $X$, but a very special
family of $p$-divisible groups: a much more rigid requirement.
See Remark \ref{john-wil} for a more on this point.

Nonetheless, the original problem has its allure and its motivation
in stable homotopy theory, and it's worth remembering this.
One flat map to $\cM_\fg$ is the   identity map
$\cM_\fg \to \cM_\fg$ itself, and we could ask whether the realization problem
can be solved for all of $\cM_\fg$.
Put aside, for the moment, the fact that $\cM_\fg$ is not an
algebraic stack, let alone a Deligne-Mumford stack. 

If the realization problem could be solved, we would have an equivalence
$$
S^0 \mathop{\longr}^{\simeq} \mathop{\holim}_{\cM_\fg} \cO_\fg^\Top
$$
where $S^0$ is the stable sphere and the homotopy limit is over  the category
of flat morphisms $\Spec(R) \to \cM_\fg$. We'd also get a descent spectral sequence
$$
H^s(\cM_\fg,\omega^{\otimes t}) \Longrightarrow \pi_{2t-s}S^0.
$$
By considering the \v Cech complex of the cover $\Spec(L) \to \cM_\fg$,
where $L$ is the Lazard ring, we have
$$
\Ext_{MU_\ast MU}^s(\Sigma^{2t}MU_\ast,MU_\ast) \cong
H^s(\cM_\fg,\omega^{\otimes t})
$$
and we would have a derived algebraic geometry version of  the Adams-Novikov
Spectral Sequence.

Alas, $\cO_\fg^{\Top}$ probably doesn't exist: we would have to use the
$fpqc$-topology on $\cM_\fg$ and, as I mentioned above, there are too many flat maps.
Nonetheless, stable homotopy theory often behaves as if $\cO_\fg^\Top$
exists. For example, the Hopkins-Ravenel chromatic convergence theorem
\cite{RavOrange} and various chromatic fracture squares are predicted
exactly by the geometry of $\cM_\fg$. A useful comparison chart is
in section 2 of \cite{HG} and much of the algebra is expanded in the later
sections of \cite{modfg}

There are four sections below. The first two sections are background on formal groups
and $\cM_\fg$. The third section makes a precise statement of the realization
problem and discusses the Hopkins-Miller theorem. This theorem states that the
realization problem has a positive answer for the moduli stack of generalized
elliptic curves. The final section discusses $p$-divisible groups and Lurie's
result.

This paper is a fleshed out version of a talk given at the conference
``New Topological Contexts for Galois Theory and Algebraic Geometry''
at the  Banff International Research Station in March of 2008. The conference
was organized by Andrew Baker and Birgit Richter. I would like to thank
the referee and the editors for careful proofreading of this paper; any
mistakes which remain are mine.

\section{Cohomology theories and formal groups}

Let's begin with a discussion of how formal groups and their invariant differentials
arise in stable homotopy theory. The following is a slight generalization of
the usual notion of a 2-periodic cohomology theory. 

\begin{defn}\label{2-period} Let $E^\ast(-)$ be a cohomology theory. Then $E^\ast$
is {\bf 2-periodic} if
\begin{enumerate}

\item the functor $X \mapsto E^\ast(X)$ is a functor to graded commutative rings;

\item for all integers $k$, $E^{2k+1}=E^{2k+1}(\ast) = 0$;

\item $E^2$ is a projective module of rank $1$ over $E^0$; and

\item for all integers $k$, the cup product map $(E^2)^{\otimes k} \to E^{2k}$
is an isomorphism.
\end{enumerate}
\end{defn}

Note that  $E^2$ is an invertible module over $E^0$ and $E^{-2}$ is the dual
module.
If $E^2$ is actually free, then so is $E_2 = E^{-2}$ and a choice of generator
$u \in E_2$ defines an isomorphism $E_0[u^{\pm 1}] \cong E_\ast$. (The
shift from $E^2$ to $E_2$ will be explained in a moment.) This happens in many important examples -- complex $K$-theory is primordial. However,
there are elliptic cohomology theories for which $E_2$ does not have a global
generator, so we insist on this generality. 

From such cohomology theories we automatically get a formal group. Recall
that if $R$ is a ring and $I \subseteq R$ is an ideal, then the {\it formal spectrum}
$$
\Spf(R,I) = \Spf(R)
$$
is the functor which assigns to each commutative ring $A$ the set of homomorphisms
$$
f:R \longr A
$$
so that $f(I)$ is nilpotent. We have an isomorphism of functors
$$
\colim \Spec(R/I^n) = \Spf(R).
$$
In many cases, $I$ is understood and dropped from the notation. Also, $R$ and
the $I$-adic completion of $R$ have the same formal spectrum, so we usually
assume $R$ is complete. As a simple example, for $\Spf(\ZZ[[x]])$ we have
 $I=(x)$ and this functor assigns to $A$ the nilpotent elements of $A$.

If $E^\ast$ is a $2$-periodic homology theory, then $E^0\CP^\infty$ is
complete with respect to the augmentation ideal 
$$
I(e) \defeq \tilde{E}^0 \CP^\infty = \mathrm{Ker}\{\ E^0 \CP^\infty \to E^0(\ast)\ \}
$$
and,  using the $H$-space structure on $\CP^\infty$, we get a commutative
group object in formal schemes
$$
G_E = \Spf(E^0 \CP^\infty).
$$
This formal group is smooth and one-dimensional in the following sense. Define
the $E_0$-module $\omega_G$ by
$$
\omega_G = I(e)/I(e)^2 \cong \tilde{E}^0S^2 \cong E_2.
$$
This module is locally free of rank 1, hence projective, and any choice of
splitting of $I(e) \to \omega_G$ defines an homomorphism out of the
symmetric algebra
$$
S_{E_0}(\omega_G) \longr E^0\CP^\infty.
$$
which becomes an isomorphism after completion.
For example, if $E_2$ is actually free
we get a non-canonical isomorphism
$$
E^0\CP^\infty \cong E^0[[x]].
$$
Such an $x$ is called a {\it coordinate}.

There is also a natural correspondence between morphisms of
cohomology theories and morphisms of formal groups.

Let $\psi:D^\ast(X) \to E^\ast(X)$ be a natural ring operation between
two 2-periodic cohomology theories. By evaluating at $X = \ast$ we
obtain a ring homomorphism $f:D^0 \to E^0$ and by evaluating at
$\CP^\infty$ we obtain a homomorphism of formal groups
$$
\phi:G_E \longr f^\ast G_D.
$$
The following result can be found in \cite{Kashiwabara} and \cite{BTurner}.
The notion of Landweber exactness is taken up below.

\begin{prop}\label{ops-homs}Let $D^\ast$ and $G^\ast$ be two Landweber
exact 2-periodic cohomology theories. Then the assignment
$$
\psi \mapsto (f,\phi)
$$
induces a one-to-one correspondence between ring operations
$$
\psi:D^\ast(-) \to E^\ast(-)
$$
and homomorphisms of pairs
$$
(f,\phi):(D^0,G_D) \to (E^0,G_E).
$$
Furthermore, $\psi$ is a stable operation if and only if $\phi$ is
an isomorphism.
\end{prop}

\begin{rem}[{\bf Formal group laws}]\label{why-fgs-not-fgls}The standard literature on chromatic
homotopy theory, such as \cite{Adams} and \cite{Rav}, emphasizes formal
group laws. If $E^\ast(-)$ is a two-periodic theory with a coordinate,
then the group multiplication
$$
G_E \times G_E = \Spf(E^0(\CP^\infty \times \CP^\infty)) \to
\Spf(E^0 \CP^\infty) = G_E
$$
determines and is determined by a power series
$$
x+_Fy = F(x,y) \in E^0[[x,y]] \cong \Spf(E^0(\CP^\infty \times \CP^\infty).
$$
This power series is a $1$-dimensional formal group law. With $2$-periodic
theories, we can insist that the formal group law be in degree zero. For 
complex oriented cohomology theories in general, the gradings become important.

Homomorphisms can also be described by power series. If $G_1$ and 
$G_2$ are two formal groups with coordinates over a base ring $R$,
then a homomorphism of formal groups $\phi:G_1 \to G_2$ is determined
by a power series $\phi(x) \in R[[x]]$ so that
$$
\phi(x+_{F_1}y) = \phi(x) +_{F_2} \phi(y)
$$
where $F_1$ and $F_2$ are the associated formal groups. The homomorphism
$\phi$ is an isomorphism if $\phi'(0)$ is a unit.
\end{rem}

\begin{rem}[{\bf Invariant differentials}]\label{invariant differential}We have defined the module
$\omega_{G_E}$ as the conormal module of the embedding
$$
e:\Spec(E^0) \to \Spf(E^0\CP^\infty) = G_E
$$
defined by the basepoint.  This definition extends to any formal group over
a base ring $R$. At first glance, this doesn't look very invariant
or very differential. We address these points.

First, $\omega_G$ has the following invariance property. If $\phi:G_1 \to G_2$
is a homomorphism of formal groups over a ring $R$, then we get an induced map
$$
d\phi:\omega_{G_2} \to \omega_{G_1}
$$
described locally as follows. If a formal group has a coordinate $x$, then
$I(e) \subseteq R[[x]]$ is the ideal of power series with $f(0)=0$ and
any element of $\omega_G$ can be written
$$
f(x)+I(e)^2 = f'(0)x + I(e)^2.
$$
Then, writing $\phi:G_1 \to G_2$ as a power series we have
\begin{equation}\label{inv1}
d\phi(f(x)+I(e)^2) = f(\phi(x))+I(e)^2 = \phi'(0)f'(0) + I(e)^2.
\end{equation}
Thus $d\phi$ is multiplication by $\phi'(0)$.

Second, while the last formula looks slightly differential, but we can do better:
$\omega_G$
is naturally isomorphic to the module of invariant differentials on $G$. 
This can be defined as follows. Let $\Omega_{G}$ denote the module
of continuous differentials on $G$; for example, if $G$ has a coordinate
$x$, then there is an isomorphism
$$
\Omega_{G} \cong R[[x]]dx.
$$
There are then three maps
$$
dp_1,dm,dp_2:\Omega_{G} \to \Omega_{G\times G}
$$
induced by the two projections and multiplication. A differential $\eta$ is
invariant if
$$
dm(\eta) = dp_1(\eta) + dp_2(\eta).
$$
Invariant differentials form an $E^0$-module; call this module $\bar{\omega}_{G_E}$
for the moment.

If $G$ has a coordinate $x$, then $\bar{\omega}_{G}$  is the free $R$-module generated
by the {\it canonical} invariant differential
$$
\eta_G = \frac{dx}{F_y(x,0)}
$$
where $F_y(x,y)$ is the partial derivative of the associated formal group
law. It is an exercise to calculate that if $\phi:G_1 \to G_2$ is a homomorphism
of formal groups with coordinate, then
$d\phi:\bar{\omega}_{G_2} \to \bar{\omega}_{G_1}$
is determined by
\begin{equation}\label{inv2}
d\phi(\eta_{G_2}) = \phi'(0)\eta_{G_1}.
\end{equation}

Finally, when $G$ has a coordinate, then we get an isomorphism
\begin{align*}
\omega_{G} = I(e)/I(e)^2 &\to \bar{\omega}_{G}\\
f(x) + I(e)^2&\mapsto f'(0)\eta_x.
\end{align*}
This isomorphisms is natural in homomorphisms (by Equations \ref{inv1}
and \ref{inv2}). In particular, it doesn't depend on
the choice of coordinate and thus
extends to the more general case as well. Because of this we drop
the notation $\bar{\omega}_{G}$
\end{rem}

\section{The moduli stack of formal groups}

Let $\cM_\fg$ be the moduli stack of formal groups: this is the algebraic geometric
object which classifies all smooth 1-dimensional formal groups and their
isomorphisms. Thus, if $R$ is a commutative ring, the morphisms
$$
G:\Spec(R) \longr \cM_\fg
$$
are in one-to-one correspondence with formal groups $G$ over $R$.
Furthermore, the 2-commutative diagrams
$$
\xymatrix@R=10pt{
\Spec(S) \ar[dr]^H \ar[dd]_f\\
&\cM_\fg\\
\Spec(R)\ar[ur]_G 
}
$$
correspond to pairs $(f:R\to S,\xymatrix{\phi:H \rto^{\cong}& f^\ast G})$.

\begin{rem}\label{stack-scheme-functor} Recall that schemes can
be defined in at least two equivalent ways. First, schemes are
defined as locally ringed spaces $(X,\cO_X)$ which have
an open cover, as locally ringed spaces, by affine schemes. This
is the point of view of Grothendieck \cite{Hart}.
Equivalently, schemes can be defined as functors from rings to 
sets which are sheaves in the Zariski topology and have an
open cover, as functors, by functors of the form
$$
A \mapsto \mathbf{Rings}(R,A).
$$
If $X$ is a scheme, in the first sense, we get a scheme in the second sense
by defining $X(R)$ to the set of all morphisms of schemes
$$
\Spec(R) \to X.
$$
This is the point of view of Demazure and Gabriel  \cite{DG}.
It is the second definition that generalizes well. A stack is then
a sheaf of groupoids on commutative rings satisfying an  
effective descent condition \cite{Laumon}\S 3. For example, $\cM_\fg$ assigns to 
each ring $R$ the groupoid of smooth one-dimensional formal
groups over $\Spec(R)$.\footnote{As in \cite{Laumon} \S 2,
we should really
speak of categories fibered in groupoids, rather than sheaves of
groupoids -- for $f^\ast g^\ast G$ is only isomorphic to $(gf)^\ast G$. However,
there are standard ways to pass between the two notions, so I will ignore the
difference.}
\end{rem}

\begin{rem}\label{covers}A scheme is more than a sheaf, of course, in
that it must have an open cover by affine schemes. Similarly, we have
{\it algebraic stacks}, which have a suitable cover
by schemes. Here is a short explanation.

A morphism $\cM \to \cN$ of stacks is {\it representable}
if for all morphisms $X \to \cN$ with $X$ a scheme, the $2$-category
pull-back (or homotopy pull-back)
$$
X \times_\cN \cM
$$
is equivalent to a scheme. A representable morphism then has algebraic
property $P$ (flat, smooth, surjective, \'etale, etc.) if all the resulting morphisms
$$
X \times_\cN \cM \to X
$$
have that property.

A stack $\cM$ is then called algebraic\footnote{The notion defined here is stronger than
what is usually called an algebraic (or Artin) stack, which requires a cover only by an algebraic
space. Algebraic spaces are sheaves which themselves have an appropriate cover by a scheme. Details are in
\cite{Laumon}.} if
\begin{enumerate}

\item every morphism $Y \to \cM$ with $Y$ a scheme is representable; and

\item there is a smooth surjective map $q:X \to \cM$ with $X$ a scheme. 
\end{enumerate}

The morphism $q$ is called a presentation. Note that an algebraic
stack may have many presentations;
indeed, flexibility in the choice of presentations leads to 
interesting theorems. See \cite{Neumann}, for an example of this
phenomenon. If the presentation can be
chosen to be \'etale, we have a {\it Deligne-Mumford stack}.
\end{rem}

\begin{rem}\label{not-alg-stack}
The stack $\cM_\fg$ is not algebraic, in this sense, as it only has
a flat presentation, not a smooth presentation. If we define $\fgl$ to
be the functor which assigns to each ring $R$ the set of formal group laws
over $R$, then Lazard's theorem \cite{Laz} says that 
$$
\fgl = \Spec(L)
$$
where $L$ is (non-canonically) isomorphic to $\ZZ[t_1,t_2,\ldots]$.
The map
$$
\fgl \longr \cM_\fg
$$
which assigns a formal group law to its underlying formal group is flat
and surjective, but not smooth since it's not finitely presented.
This difficulty can be surmounted in two ways: enlarge the notion of an
algebraic
stack to include flat presentations or note that $\cM_\fg$ can be
written as the 2-category inverse limit of a tower of the algebraic
stacks of ``buds'' of formal groups and is, thus, pro-algebraic.
See \cite{Smith}.
\end{rem}

\begin{rem}\label{mod-sheaves} A sheaf in the $fpqc$-topology
on $\cM_\fg$ is a functor $\cF$ on the category of affine
schemes over $\cM_\fg$ which satisfies faithfully flat
descent. Thus, for each formal group $G$
over $R$ we get a set (or ring, or module, etc.) 
$$
\cF(R,G) = \cF(G:\Spec(R) \to \cM_\fg)
$$
and for each $2$-commuting diagram
$$
\xymatrix@R=10pt{
\Spec(S) \ar[dr]^H \ar[dd]_f\\
&\cM_\fg\\
\Spec(R)\ar[ur]_G 
}
$$
a restriction map $\cF(S,H) \to \cF(R,G)$. This must be a sheaf in 
the sense that if $q:S \to R$ is faithfully flat, then there is an equalizer
diagram
$$
\xymatrix{
\cF(R,G) \rto & \cF(S,q^\ast G) \ar@<.5ex>[r]  \ar@<-.5ex>[r]
& \cF(S \otimes_R S,p^\ast G)
}
$$
where I have written $p$ for the inclusion $R \to S \otimes_R S$.

For example, define the structure sheaf $\cO_\fg$
to be the functor on affine schemes over $\cM_\fg$ with
$$
\cO_\fg(R,G) = \cO_\fg(G:\Spec(R) \to \cM_\fg) = R.
$$
More generally, we consider module sheaves $\cF$ over $\cO_\fg$. Such
a sheaf is called {\it quasi-coherent} if, for each $2$-commutative
diagram, the restriction map $\cF(R,G) \to \cF(S,H)$ extends to
an isomorphism
$$
S \otimes_R \cF(R,G) \cong \cF(S,H).
$$
This isomorphism can be very non-trivial, as it depends on the 
choice of isomorphism $\phi:H \to f^\ast G$ which makes the diagram
$2$-commute.

A fundamental example of a quasi-coherent sheaf is the sheaf 
of invariant differentials $\omega$ with
$$
\omega(R,G) = \omega_G
$$
the invariant differentials on $G$. This is locally free of rank $1$ and hence
all powers $\omega^{\otimes n}$, $n \in \ZZ$, are also quasi-coherent sheaves.
The effect of the choice of isomorphism in the $2$-commuting diagram on the
transition maps for $\omega^{\otimes n}$ is displayed in Equations
\ref{inv1} and \ref{inv2}.
\end{rem}

\begin{rem}Consider the $2$-category pull-back
$$
\fgl \times_{\cM_\fg} \fgl
$$
where $\fgl =\Spec(L)$ is the functor of formal group laws. The pull-back
is the functor which assigns to each commutative ring the set of triples
$(F_1,F_2,\phi)$ where $F_1$ and $F_2$ are formal group laws and 
$\phi$ is an isomorphism of their underlying formal groups. Given that
these formal groups have a chosen coordinate, the isomorphism $\phi$ can
be expressed as an invertible power series $\phi(x) = a_0 x + a_1 x^{2} + \cdots$.
Thus the pull-back is the affine scheme on the ring
$$
W = L[a_0^{\pm 1},a_1,\ldots].
$$
The pair $(L,W)$ forms a {\it Hopf algebroid}; that is, a groupoid in affine
schemes.
Furthermore, the category of quasi-coherent sheaves on $\cM_\fg$ is
equivalent to the category of $(L,W)$-comodules. 

To get a functor in one direction, let $\cF$ be a quasi-coherent sheaf.
Then
$$
M = \cF(\fgl \to \cM_\fg)
$$
is an $L$-module. One of the two projections $\fgl \times_{\cM_\fg} \fgl \to \fgl$
shows
$$
 \cF(\fgl \times_{\cM_\fg} \fgl \to \cM_\fg) \cong W \otimes_L M
$$
and the other projection supplies the comodule structure map.
I will say something about how you pass from comodules to sheaves
at the beginning of the next section.

It is here we see the flexibility of choosing the presentation. For example, 
if we work localized at some prime $p$ and consider the stack
$$
\Spec(\ZZ_{(p)}) \times \cM_\fg \defeq \ZZ_{(p)} \otimes \cM_\fg
$$
then we could use the scheme of $p$-typical formal group laws as our
cover and obtain a different category of comodules closely related (up to
issues of grading) to $(BP_\ast,BP_\ast BP)$-comodules. This would then
be equivalent to the category of $\ZZ_{(p)} \otimes (L,W)$-comodules
as both would be equivalent to quasi-coherent sheaves on
$\ZZ_{(p)} \otimes \cM_\fg$.
\end{rem}

\begin{rem}[{\bf The height filtration}]\label{height} If $G$ is
a formal group over a field $\FF$ we can always choose a coordinate.
If $\FF$ is of characteristic $p$ for some prime $p$, then the
homomorphism $p:G \to G$ can be written
$$
p(x) = ux^{p^n} + \cdots
$$
for some $u \ne 0$ and $1 \leq n \leq \infty$. (If $n=\infty$, $p = 0:G \to G$.)
The number $n$ is an isomorphism invariant and, if $\FF$ is separably
closed, a complete invariant, by the theorem of Lazard \cite{Laz}.
This notion of {\it height} can be extended to formal groups over an
arbitrary $\FF_p$-algebra or even to formal groups over schemes
over $\FF_p$, but some care is needed if $G$ does not have a coordinate.

Consider a formal group $G$ over an $\FF_p$-algebra $R$.
If we let $f:R \to R$ be the Frobenius, we get a new formal group
$G^{(p)} = f^\ast G$. We then have a diagram
$$
\xymatrix{
G \ar[dr] \rto_F \ar@/^1pc/[rr]^f  &G^{(p)}\dto \rto & G \dto\\
&\Spec(R) \rto_f &\Spec(R)
}
$$
where the square is a pull-back. The homomorphism $F$ is the
{\it relative Frobenius}. We know that if
$$
\phi:G \to H
$$
is a homomorphism of formal groups over $R$
for which $d\phi = 0:\omega_H \to \omega_G$, there is then a factoring
$$
\xymatrix{
G \rto_F  \ar@/^1pc/[rr]^\phi& G^{(p)} \rto_\psi &H
}
$$
Then we can test $d\psi$ to see if we can factor further.\footnote{See \cite{modfg}\S 5 for a proof
of these facts in this language.}

For example, let $\phi = p:G \to G$ be $p$th power map. Then we obtain
a factoring 
$$
\xymatrix{
G \rto_F  \ar@/^1pc/[rr]^p& G^{(p)} \rto_{V_1} &G
}
$$
This yields an element
$$
dV_1 \in \Hom(\omega_G,\omega_{G^{(p)}})
$$
and we can factor further if $dV_1=0$. Since $G$ is of dimension 1,
$\omega_{G^{(p)}}\cong \omega_G^{\otimes p}$; since $\omega_G$ is invertible,
$$
 \Hom(\omega_G,\omega_G^{\otimes p}) = \Hom(R,\omega_G^{\otimes p-1}).
$$
Thus $dV$ defines an element\footnote{Over a base scheme
which was not affine, $v_1(G)$ is a global section.} $v_1(G)$ of $\omega_G^{\otimes p-1}$. If $v_1(G) = 0$, then we obtain
a further factorization
$$
\xymatrix{
G \rto_F  \ar@/^1pc/[rrr]^p& G^{(p)} \rto_{F^{(p)}} & G^{(p^2)} \rto_{V_2} &G
}
$$
and an element $v_2(G) \in \omega_{G}^{\otimes p^2-1}$. This can be continued
to define elements $v_n(G) \in \omega_G^{\otimes p^n-1}$ and $G$ has height
at least $n$ if $v_1(G) = \cdots = v_{n-1}(G) = 0$. We say $G$ has height
{\it exactly} $n$ if $v_n(G) \in \omega_G^{\otimes p^n-1}$ is a generator.

The elements $v_n(G)$ defined in this way are isomorphism invariants. For example
if $G$ is a formal group over an $\FF_p$ algebra $R$, then 
$$
p(x) = u_1x^p + \cdots.
$$
The element $u_1$ is not an isomorphism invariant, but if $\eta= dx/F_y(x,0)$ is the
standard invariant differential, then
$$
v_1(G) = u_1 \eta^{\otimes p-1} \in \omega_G^{\otimes p-1}
$$
is an invariant. 

Because of this invariance property, the assignment $G \mapsto v_1(G)$
defines a global section $v_1$ of the sheaf $\omega^{\otimes p-1}$ on
the closed substack
$$
\FF_p \otimes \cM_\fg \defeq \cM(1) \subseteq \cM_\fg
$$
In this way we obtain a sequence of closed substacks
$$
\dots \subseteq \cM(n+1) \subseteq \cM(n) \subseteq
\cdots \subseteq \cM(1) \subseteq \cM_\fg
$$
where $\cM(n+1) \subseteq \cM(n)$ is defined by the vanishing of the
global section $v_n$ of $\omega^{\otimes p^n-1}$. Thus $\cM(n)$ classifies
formal groups of height at least $n$. The relative open
$$
\cH(n) = \cM(n) - \cM(n+1)
$$
classifies formal groups of height exactly $n$. Lazard's theorem, rephrased,
says that $\cH(n)$ has a single geometric point given by a formal group
$G$ of height $n$ over any algebraically closed field $\FF$. The pair
$(\FF,G)$ has plenty of isomorphisms, however, so $\cH(n)$ is not a scheme;
indeed, in the language of \cite{Laumon} it is a {\it neutral gerbe}.
See \cite{Smith}.

One way to think of neutral gerbes is as the {\it classifying stack}
which assigns to any commutative
ring $R$ the groupoid of torsors for some group scheme. In this case, we can
take the following group scheme. Let $\Ga_n$ be some fixed formal group
defined over the finite field $\FF_p$ and let $\Aut(\Ga_n)$ be the
group scheme which assigns to each $\FF_p$ algebra $R$ the group
of automorphisms of $i^\ast \Ga_n$ over $R$. Here I write $i:\FF_p \to R$
for the inclusion. Then $\cH(n)$ is the classifying stack for $\Aut(\Ga_n)$.

This group scheme is quite familiar to homotopy theorists. To be specific, let
$\Ga_n$ be the Honda formal group over $\FF_p$; this is the $p$-typical formal group
law with $p(x) = x^{p^n}$. Then $\Aut(\Ga_n)$
is the affine group scheme obtained from the Morava stabilizer algebra
(See \S 6.2 of \cite{Rav}) and the group of
$\FF_{p^n}$-points of $\Aut(\FF_p,\Ga_n)$ is the Morava
stabilizer group $S_n$.
By definition, $S_n$ is the automorphisms of $\Ga_n$ over $\FF_{p^n}$.
\end{rem}

\begin{rem}[{\bf Landweber's criterion for flatness}]\label{LEFT} We will be
concerned with morphisms $\cN \to \cM_\fg$ of stacks which are representable
and flat. We defined this notion above, but the Landweber exact functor
theorem  uses the height filtration to give an easily checked criterion
for flatness. This we now state, first giving a global formulation, then giving
a way to check this locally.

We begin by noting that the closed inclusion $j:\cM(n+1) \to \cM(n)$
is actually an effective Cartier divisor. This means the following. Let $\cO(n)$ be
the structure sheaf of $\cM(n)$. Then the 
global section $v_n \in H^0(\cM(n),\omega^{\otimes p^n-1})$
defines an {\it injection} of sheaves
$$
\xymatrix{
0 \to \cO(n) \rto^-{v_n} \rto &\omega^{\otimes p^n-1} }.
$$
This yields a short exact sequence
$$
\xymatrix{
0 \to \omega^{\otimes -(p^n-1)} \rto^-{v_n} &\rto \cO(n)\rto & j_\ast \cO(n+1) \to 0.}
$$
This identifies $\omega^{\otimes -(p^n-1)}$ with the ideal defining 
the closed inclusion $\cM(n+1)$ in $\cM(n)$.

Now let $f:\cN \to \cM_\fg$ be a representable morphism of stacks and let
$$
\cN(n) = \cM(n) \times_{\cM_\fg} \cN \subseteq \cN.
$$
Then $\cN(n+1) \subseteq \cN(n)$ remains a closed inclusion and
if $f$ is flat, it remains an effective Cartier divisor; that is,
$$
\xymatrix{
0 \to \cO_{\cN(n)} \rto^-{v_n} & \omega^{\otimes p^n-1}
}
$$
remains an injection. Landweber's theorem \cite{Land} now
says that this is sufficient.\footnote{For proofs in the language
of stacks see \cite{LEFT}, \cite{modfg}, \cite{Neumann}, and \cite{Hollander}.
The first two of these were inspired directly by Mike Hopkins, the third had
 input from Mark Behrens.} That is, suppose that for all primes $p$
and all integers $n$, the morphism
$$
v_n:\cO_{\cN(n)} \to \omega^{\otimes p^n-1}
$$
is an injection.
Then $f: \cN \to \cM$ is flat. For proofs in the language
of stacks see \cite{LEFT}, \cite{Neumann}, and \cite{Hollander}.
The first of these (which has an extra hypothesis) was inspired directly
by Mike Hopkins, the second had input from Mark Behrens.

Locally this can be unwound as follows. Let $\Spec(R) \to \cM_\fg$
be a formal group with a coordinate $x$. Define $u_0 = p$
and recursively define elements $u_n$ by
$$
p(x) = u_nx^{p^n} + \cdots
$$
modulo $(p,u_1,\ldots,u_{n-1})$. Then we can rewrite the equations
above as saying that, for all primes $p$ and all $n$, the
multiplication
$$
u_n:R/(p,u_1,\ldots,u_{n-1}) \to R/(p,u_1,\ldots,u_{n-1})
$$
is an injection. 

Much of the proof of Landweber's result is formal, using only that the closed
inclusions $\cM(n+1) \subseteq \cM(n)$ are effective Cartier divisors.
But in the
end, one must use something about formal groups, and -- as Neil Strickland
has pointed out -- the crucial ingredient
turns out to be Lazard's uniqueness theorem in the following strengthened form.
\end{rem}

\begin{prop}\label{laz-strong}Let $G_1$ and $G_2$ be two formal
groups of the same height over an $\FF_p$ algebra $R$. Then there is a sequence of
\'etale extensions
$$
R \subseteq R_1 \subseteq R_2 \subseteq \cdots
$$
so that $G_1$ and $G_2$ become isomorphic over $R_\infty = \colim R_n$.
\end{prop}

This is what is actually proved by Lazard in \cite{Laz}.
See also  \cite{Rav}, Appendix 2.2 and
\cite{modfg}, Theorem 5.25 for this result over arbitrary base schemes.
There are additional references in \cite{modfg}. If
$R=\FF$ is field, the extension adjoins roots of certain separable
polynomials; hence if $\FF$ is separably closed, $G_1$ and $G_2$
were already isomorphic.

\section{The realization problem}

The Landweber Exact Functor Theorem was originally proved to provide
homology theories. This begins with periodic complex cobordism
$MUP_\ast$, which is obtained from ordinary complex cobordism by 
adjoining an invertible element of degree $2$:
$$
MUP_\ast X = \ZZ[u^{\pm 1}] \otimes_\ZZ MU_\ast X.
$$
The representing spectrum is $\vee_n \Sigma^{2n}MU$,  the Thom spectrum
of the universal bundle over $\ZZ \times BU$. Notice that the wedge summands
keep track of the virtual dimension over the individual components. We
have
$$
MUP_0 = L\qquad\mathrm{and}\qquad MUP_0MUP = W.
$$

Now suppose we are given a ring $R$ and a formal group $G$ with a coordinate
over $R$. The choice of coordinate defines a map of rings $L \to R$
and we can examine the functor
$$
X \mapsto R \otimes_L MUP_\ast X.
$$
Landweber's criterion guarantees that this functor yields a homology theory
$E(R,G)_\ast$. A theorem of Hovey and Strickland \cite{HS} says that there are no
phantom maps between such theories and this implies that $E(R,G)$ is
actually a homotopy commutative ring spectrum. To get the multiplication
map, for example, note that
$$
E(R,G)_0E(R,G) \cong R \otimes_L W \otimes_L R
$$
still satisfies Landweber's criterion; hence the morphism
$$
R \otimes_L W \otimes_L R \to R
$$
which classifies the identity from $G$ to itself defines  the multiplication.
This is as natural as can be:  we get a functor from
formal groups with coordinate to the stable homotopy category.

I'd next like to eliminate the reliance on the coordinate. A formal group
need not have a coordinate and, even if it does, I'd rather not choose one.

\begin{rem}[{\bf From comodules to sheaves}]\label{from-comods-to-sheaves}
Suppose we are given an $(L,W)$ comodule $M$; we'd like to produce
a quasi-coherent sheaf $\cF_M$ on $\cM_\fg$. Let $G:\Spec(R) \to \cM_\fg$
be flat. Consider the diagram with all squares pull-backs:
$$
\xymatrix{
X_1 \ar@<.5ex>[r]\ar@<-.5ex>[r]\dto & X_0 \rto^-f \dto & \Spec(R)\dto^G\\
\fgl \times_{\cM_\fg}\fgl \ar@<.5ex>[r]\ar@<-.5ex>[r]\ & \fgl \rto_-q & \cM_\fg.
}
$$
Here I have written $\fgl$ for $\Spec(L)$ and $\fgl\times_{\cM_\fg}\fgl$
for $\Spec(W)$. The scheme $X_0$ represents the functor which assigns
to each commutative $R$-algebra $A$ the set of coordinates of $G$
over $A$. If I could choose a coordinate for $G$, we'd get a non-canonical
isomorphism
$$
X_0 \cong \Spec(W \otimes_L R).
$$
Since such a choice is always possible locally, we conclude $f$ is
an affine morphism of schemes.  Since $q$ is faithfully flat, so is $f$ and
to specify $\cF_M(R,G)$ I need only specify a quasi-coherent sheaf
on $X_0$ together with descent data. This sheaf is the pull-back of the sheaf
on $\fgl$ determined by $M$; the descent data is determined by the
comodule structure and the commutative diagram.

This elaborate description is a choice-free way of naming $\cF(R,G)$. If
we can choose a coordinate for $G$, then we get an isomorphism
$$
\cF_M(R,G) \cong R \otimes_L M.
$$
\end{rem}
\def\modfg{{\cM_\fg}}

Now suppose $\Spec(R) \to \modfg$ is flat and classifies the formal group
$G$. Define a homology theory by 
$$
E(R,G)_\ast X = \cF_{MUP_\ast X}(R,G)
$$
where $\cF_M$ is the sheaf associated to the comodule $M$.
By Hovey and Strickland's result, quoted above, this is a homotopy commutative
ring theory. Furthermore, $E(R,G)_0 \cong R$, $E(R,G)_{2k+1}=0$, the
associated formal group is $G$ and $E(R,G)_{2k} \cong \omega_G^{\otimes k}$.
In this way, we get a presheaf
$$
E(-,-):\mathrm{Flat}/\modfg \longr \mathbf{Ho}(\mathrm{Spectra})
$$ 
from the category of flat maps with affine source over $\modfg$ to 
the stable homotopy category
realizing the graded structure sheaf $\cO_\ast = \{\omega^{\otimes \ast}\}$.
Here and throughout we assume $\omega$ is in degree $2$, for topological
reasons.

The realization problem asks to what extent the presheaf $E(-,-)$ can
be lifted to the category of $E_\infty$ ring spectra. 

Since the geometry of $\modfg$ is not so good -- it is not an algebraic stack,
for example (see \ref{not-alg-stack}) --  I don't suppose anyone
expects an affirmative answer to this question -- there are simply too many flat maps. (For further comments on this point, see the introduction.)
One way to cut down the class of morphisms
is to restrict attention to stacks with more structure. Over an algebraic stack,
for example, we can work with the smooth-\'etale topology (\cite{Laumon},
\S 12) and over a Deligne-Mumford stack we can work with
the \'etale topology.  Thus, I will formulate the question as follows:

\begin{realize} Let $\cM \to \cM_\fg$ be a representable and
flat morphism from an algebraic stack and let
$$
\cO_{\cM\ast} = \{\omega^{\otimes \ast}\}
$$
denote the graded structure sheaf on $\cM$ in an appropriate topology.
Is there a pre-sheaf of $E_\infty$-ring
spectra $\cO^\Top_S$ with an isomorphism of associated sheaves
$$
(\cO^\Top_S)_{\ast} \cong \cO_{\cM\ast}?
$$
If so, how unique is this? What is the homotopy type of the space of all
realizations?
\end{realize}

\begin{exam}\label{LEFT-realize}Even in this generality, the problem might
not have a general solution. For example, we could take $\cM = \Spec(R)$ 
and $G:\Spec(R) \to \cM_\fg$ to be any flat map and the Zariski topology
on $\Spec(R)$. Then a positive solution
to the realization problem would say that representing spectrum $E(R,G)$
of the resulting Landweber exact homology theory had the structure of
an $E_\infty$-ring spectrum. This is not very likely.
More on  this point
below in Remark \ref{john-wil}.
\end{exam}

There is a very important example of a positive solution of the realization
problem: the moduli stack of elliptic curves. Standard references on
elliptic curves include \cite{Silver} and \cite{KM}; the stack was introduced
in \cite{DM} and thoroughly studied in \cite{DR}.

\begin{rem}[{\bf The moduli stack of elliptic curves}]\label{mell} Let $S$
be a scheme. Then an elliptic curve over $S$
$$
\xymatrix{
C\ar@<.5ex>[r]^-q & \ar@<.5ex>[l]^-e S
}
$$
is a proper, smooth curve over $S$, with geometrically connected fibers
of genus $1$ and with a given section $e$. Such curves have a natural
structure as an abelian group scheme $S$ with $e$ as the identity section.
By taking
a formal neighborhood of $e$ in $C$ we obtain a formal group $C_e$.
There is a stack $\mell$ (called $\cM_{1,1}$ in the algebraic geometry
literature -- ``genus 1 with 1 marked point'') classifying elliptic curves;
that is, morphisms
$$
C:S \to \mell
$$
are in one-to-one correspondence with elliptic curves over $S$. This
stack was produced in $\cite{DM}$ and is one of the original examples
of a stack. The assignment $C \to C_e$ produces a morphism of stacks
$$
\mell \longr \modfg
$$
which is representable and flat and we could ask about the realization
problem for $\mell$. The sheaf $\omega$ on $\modfg$ restricts to the
sheaf on $\mell$ which assigns to each elliptic curve $C$ over $S$
the sheaf $\omega_C$ on $S$ of invariant differentials of $C$. The global
sections
$$
H^0(\mell, \omega^{\otimes t})
$$
are the modular forms of weight $t$ (and level one); they assemble into
a graded ring. From \cite{DelMod} we have an isomorphism
$$
\ZZ[c_4,c_6,\Delta^{\pm 1}]/(c_4^3 - c_6^2 = (12)^3\Delta) \cong
H^0(\mell,\omega^{\otimes \ast}).
$$
where $c_4$, $c_6$, and $\Delta$ are the standard modular
forms of weight (degrees) 4, 6, and 12 respectively.
\end{rem}

\begin{rem}[{\bf The compactification of $\mell$}]\label{DR-compact}
There is a canonical compactification
$\barell$  of the moduli stack $\mell$. One way to construct this as follows.

Locally in $S$ any elliptic
curve is a non-singular subscheme of $\PP^2$ obtained from a Weierstrass
equation
$$
y^2 + a_1xy + a_3y = x^3+ a_2x^2 + a_4x + a_6.
$$
Any such curve is called a Weierstrass curve; more generally, we define
a Weiertrass curve $C$ over a scheme $S$ to be a pointed morphism of schemes
$\xymatrix{C \ar@<.5ex>[r] & \ar@<.5ex>[l]^e S}$ which can be given
Zariski-locally by a Weierstrass equation. (The marked point of a
Weierstrass curve is the point $e=[0,1,0]$).
Although not every Weierstrass curve is an elliptic curve,
we do get an embedding 
$$
\mell \longr \mweier
$$
into a moduli stack of Weierstrass curves and the morphism
$\mell \to \modfg$ factors through this embedding:
$$
\mell \to \mweier \to \modfg.
$$
This is a consequence of the fact that, for any Weierstrass curve $C$, the marked point
$e=[0,1,0]$ is always smooth and the smooth locus on $C$ has a natural structure
as an abelian group scheme with $e$ as the identity.

The morphism $\mweier \to \modfg$ is not good, however, for two reasons:
first, it is not flat (see Rezk \cite{512}\S 20) and, second, the geometry
of $\mweier$ is not very good. For example, because the automorphism group
of the cusp curve $y^2 = x^3$ is not an \'etale group scheme, this stack
cannot be a Deligne-Mumford stack. (See \cite{Laumon}, Th\'eor\`eme 8.1.)
However, the sheaves $\omega^{\otimes t}$ yield
sheaves on $\mweier$ and a calculation from \cite{512} and \cite{Bauer}
(following Deligne \cite{DelMod}, of course) 
implies that there is an isomorphism of graded rings
$$
\ZZ[c_4,c_6,\Delta]/(c_4^3 - c_6^2 = (12)^3\Delta) \cong
H^0(\mweier,\omega^{\otimes \ast}).
$$
A Weierstrass curve $C$ is an elliptic curve  and, hence,
smooth if $\Delta(C)$ is invertible. We define
$$
\barell \longr \mweier
$$
to be the substack of curves $C$ so that the sections
$c_4^3(C)$, $c_6^2(C)$, and $\Delta(C)$ generate $\omega_C^{\otimes 12}$;
that is, in formulas, we have:
$$
(c_4^3(C),c_6^2(C),\Delta(C)) = \omega_C^{\otimes 12}.
$$
There is an inclusion $\mell \subseteq \barell$; however, we also
allow other curves -- for example, we allow curves where $c_4(C)$ 
is invertible. In effect, we allow nodal,
but not cusp, singularities. Thus
$$
y^2 = x^2(x-1)
$$
is allowed, but $y^2=x^3$ is not.\footnote{One important
generalized elliptic now allowed is the Tate curve over
$\ZZ[[q]]$ which is singular at $q=0$. See \cite{DR}, \S VII.
The morphism $\Spf(\ZZ[[q]]) \to \barell$ classifying the Tate
curve identifies $\Spf(\ZZ[[q]])$ as formal neighborhood of the
singular generalized elliptic curves. }
The resulting map
$$
\barell \longr \modfg
$$
is flat and $\barell$ has good geometry.
\end{rem}

\begin{rem}[{\bf \'Etale maps to $\barell$}]\label{etale-mell}
To get some feel for the
geometry of $\barell$, define the $j$-invariant
\begin{align*}
j:\barell &\longr \PP^1\\
C \longmapsto &[c_4^3(C),\Delta(C)].
\end{align*}
Then from \cite{DR} \S V1.1 we learn that $j$
identifies $\PP^1$ as the the ``coarse moduli stack'' of $\barell$ --
the scheme which most closely approximates the sheaf  of
isomorphism classes of generalized
elliptic curves. (See \cite{DR}, \S I.8 for precise definitions.)
Furthermore, the fiber at any $j$-value is
the classifying stack of a finite group scheme.

Note that while the geometry of the $j$-invariant is somewhat
complicated, $\barell$ is still smooth of dimension $1$. (See
\cite{DR}.) In particular, if $\Spec(R) \to \barell$ is \'etale, then
$R$ is smooth of dimension 1 over $\ZZ$. There are
classical examples of such affine morphisms; see, for example,
\cite{Igusa}. The
one emphasized in usual sources (see \cite{Silver} \S III.1) is the
Legendre curve 
$$
y^2 = x(x-1)(x-\la)
$$
over $\ZZ[1/2][\la, (\la^2 - \la)^{-1}]$. There is also the Deuring curve
(\cite {Silver} Proposition A.1.3)
$$
y^2 + 3\nu + xy = x^3
$$
over $ \ZZ[1/3][\nu,1/(\nu^3+1)]$. Both the Legendre curve and the
Deuring curve are smooth, because we've inverted the discriminant
$\Delta$. A wide example of non-smooth curves can be obtained by
base change from the curve
$$
y^3 + xy = x^3 + \tau
$$
over $\ZZ[\tau,1/(1 + 2^43^3\tau)]$. (We invert $1 + 2^43^3\tau$ to make
the singular locus of this curve exactly $\tau=0$.) An observation, which
I learned from Hopkins, is that these three curves form an affine \'etale
cover of $\barell$.
\end{rem}

Here is the famous positive answer to the realization problem. See
\cite{ICM}.

\begin{thm}[{\bf Hopkins-Miller}]\label{HM} The realization problem for
$$
\barell \longr \modfg
$$
has a solution in the \'etale topology: there is a presheaf $\cO_\Ell^\Top$ of
$E_\infty$-ring spectra realizing the graded structure sheaf $\cO_{\Ell\ast}$. 
The space of all realizations is path connected.
\end{thm}
\def\TMF{{{\mathbf{tmf}}}}

If we define $\TMF$ to be the homotopy global sections
$$
\TMF \defeq \mathop{\holim}_{\barell} \cO_\Ell^\Top
$$
where the homotopy limit is over all \'etale morphisms
$\Spec(R) \to \barell$.
There is a descent spectral sequence
\begin{equation}\label{descent}
H^s(\barell,\omega^{\otimes t}) \Longrightarrow \pi_{2t-s}\TMF
\end{equation}
and modular forms are, by definition,
$$
H^0(\barell,\omega^{\otimes \ast}) \cong H^0(\mweier,\omega^{\otimes \ast}) \cong
\ZZ[c_4,c_6,\Delta]/(c_4^3-c_6^2=(12)^3\Delta).
$$
Hence ``topological modular forms''. The question of which modular forms are
homotopy classes is quite interesting. For example, $c_6$ and the discriminant
$\Delta$ are not; however $2c_6$ and $24\Delta$ are. See \cite{ICM}.

The calculation of $\pi_\ast \tmf$ has been made completely. While not
yet explicitly in print, it can be easily deduced from \cite{512} and \cite{Bauer} --
both of which follow \cite{HopMah}. There is a curious feature of the answer:
while the $E_2$ term does not display any obvious duality, the homotopy groups
of $\tmf$ have a very strong duality very similar to Serre duality for 
projective schemes. I know of no good explanation for this -- the differentials
and extensions in the spectral sequence conspire in an almost miraculous
fashion to give the result -- but there must be one in derived algebraic
geometry. Compare also Mahowald-Rezk duality \cite{MahRez}.

\begin{warn}  Note that topological modular forms has often been
defined to be the zero-connected
cover of what I've called $\TMF$. However, to even decide if this makes
sense, you need to calculate $\pi_\ast\tmf$ and notice that the resulting
answer takes a very special form. 
\end{warn}

\begin{exam}[{\bf Topological Automorphic Forms}]\label{taf} Work
of Mark Behrens and Tyler Lawson \cite{BL} solve the realization problem for
certain Shimura varieties, which are moduli stacks of highly structured
abelian varieties. The extra structure is needed to get formal groups
of higher heights. The problem is that the only abelian group schemes
of dimension one are the additive group $\GG_a$, the multiplicative
group $\GG_m$, and elliptic curves; from these we only get 
formal groups of height $\infty$, $1$, and $2$. To get formal groups
of height greater than $2$, one must use higher dimensional abelian
group schemes $A$, but then one must add enough structure so that 
the formal completion $A_e$ of $A$ at the identity splits off a natural
summand of dimension one. It takes a while to define such objects -- so I
won't do it here -- but it turns out they've been heavily studied in number theory.
See, for example, \cite{Kottwitz}.
\end{exam}

\begin{rem}[{\bf The role of $E_\infty$-ring spectra}]\label{why-einfty} Why do I 
(following my betters, notably Mike Hopkins) insist on highly
structured ring spectra in the realization problem? There are
two reasons.
\begin{enumerate}
\item {\bf (Practical)} Asking for $E_\infty$-ring spectra allows for algebraic geometry 
(i.e., ring theoretic) input into the constructions; and

\item {\bf (Aesthetic)} The stack $\mell$ with its $E_\infty$ structure sheaf becomes 
a central exhibit in the world of derived algebraic geometry: we learn something
inherently new about elliptic curves.
\end{enumerate}
\end{rem}

\section{Lurie's theorem and p-divisible groups}

The formal group of an elliptic curve is part of a richer and more rigid structure.
At this point we pick a prime $p$ and work over $\Spf(\ZZ_p)$; that is,
$p$ is implicitly nilpotent in all our rings. This has the implication that
we will we working in the $p$-complete stable category.

\begin{defn}\label{pdiv} Let $R$ be a ring and $G$ a sheaf of abelian groups
on $R$-algebras. Then $G$ is a {\bf $p$-divisible group} of {\bf height $n$}
if 
\begin{enumerate}

\item $p^k:G \to G$ is surjective for all $k$;

\item $G(p^k) = \mathrm{Ker}(p^k:G \to G)$ is a finite and flat group
scheme over $R$ of rank $p^{kn}$;

\item $\colim G(p^k)\cong G$.
\end{enumerate}
\end{defn}

\def\et{{{\mathrm{{et}}}}}

\begin{rem}\label{p-div-rem} 1.) If $G$ is a $p$-divisible group, then
completion at $e \in G$ gives an abelian formal group $G_\for \subseteq G$,
not necessarily of dimension $1$. The quotient $G/G_\for$ is 
\'etale over $R$; thus we get a natural short exact sequence
$$
0 \to G_\for \to G \to G_\et \to 0.
$$
This is split over fields, but not in general.

2.) If $C$ is a smooth elliptic curve, then $C(p^\infty) = \colim C(p^n)$ is
$p$-divisible of height $2$ with formal part of dimension 1.

3.) Formal groups need not be $p$-divisible groups as there is no reason to suppose
$$
\colim G(p^k) \cong G
$$
over a ring which is not local and complete. Nor can one assume that a formal
group is a sub-group a $p$-divisible group. 

4.) (Rigidity) If $G$ is a $p$-divisible group over a scheme $S$, the function
which assigns to each geometric point $x$ of $S$ the height of the
fiber $G_x$ of $G$ at $x$ is constant. This is not true of formal groups, as
the example of elliptic curves shows. Indeed, if $G$ is $p$-divisible
of height $n$ with $G_\for$ of dimension 1, then the height of $G_\for$
can be any integer between $1$ and $n$.

For a simple example of this phenomenon, take $n>1$ and let 
$$
S = \Spec(\ZZ/(p^n)[u_1])
$$
and $G$ the formal group obtained from the  $p$-typical formal group
law $F$ with $p$-series
$$
p_F(x) = px +_F u_1x^p  +_F x^{p^2}.
$$
Then if $x$ is the point given by the maximal ideal $(p,u_1)$, $G_x$
has height $2$;  however, if $x$ is point given the ideal $(p,u_1-1)$,
then $G_x$ has height $1$. This example is closely related to the
Johnson-Wilson theory $E(2)_\ast$. It is not at all clear that this
formal group has anything to do with a $p$-divisible group.
See, however, cite{HL}, where the authors do have some success at the
prime $3$.
\end{rem}

\begin{exam}[{\bf $p$-divisible groups and localization}] The following example,
which I learned from Charles Rezk, shows that $p$-divisible groups arise
naturally in homotopy theory.

Let $E=E_n$
be a Morava $E$-theory; this is a $2$-periodic theory with a non-canonical
isomorphism
$$
E_0 = W(\FF_{p^t})[[u_1,\dots,u_{n-1}]]
$$
and whose formal group is  a universal deformation of a height $n$-formal
group.  (See Examples \ref{ST} and \ref{Def} below for more on deformations.)
Since $E_0$ is complete
$$
G = G_{E} = \Spf(E^0\CP^\infty)
$$
is a $p$-divisible group of height $n$. Indeed,
$$
\map(\CP^\infty,E) \simeq \map(\colim BC_{p^k},E)
\simeq \lim \map(BC_{p^k},E)
$$
and, by applying $\pi_0$ we get
$$
\Spf(E^0\CP^\infty) \cong \colim \Spec(E^0BC_{p^k}) = \colim G(p^k).
$$
When we apply the localization functor $L_{n-1} = L_{E(n-1)}$, we
have
$$
\xymatrix@C=15pt{
\map(\CP^\infty,L_{n-1}E) & \lto L_{n-1}\map(\CP^\infty,E) \rto & 
L_{n-1}\map(BC_{p^k},E)
}
$$
yielding
$$
G_{{}_{L_{n-1}E}}\ \longr \colim \Spec(\pi_0L_{n-1}\map(BC_{p^k},E))
$$
as the inclusion of the formal part of a $p$-divisible group. This map
is not an isomorphism; indeed,  the rank of
$$
G_{{}_{L_{n-1}E}}(p)
$$
over $\pi_0L_{n-1}E$ is $p^{n-1}$ while the rank of 
$$
\Spec(\pi_0L_{n-1}\map(BC_{p},E))
$$
is $p^n$. This last group scheme is the $p$-torsion
in the $p$-divisible group.

Dealing
with examples such as this is one of the deeper technical aspects of the
original Hopkins-Miller proof of the existence of $\TMF$.
\end{exam}

\begin{defn}\label{mod-p-div} Let $\cM_p(n)$ be the moduli stack
of $p$-divisible groups of
height $n$ and with $\mathrm{dim}\ G_\for = 1$.
\end{defn}

\begin{rem}\label{p-div-rems}The stack $\cM_p(n)$ is not an algebraic stack,
but rather pro-algebraic. This can be deduced from the material
in the first chapter of \cite{Messing}.
\end{rem}

\begin{rem}\label{tomodfg}There is a morphism of stacks
\begin{align*}
\cM_p(n) &\longr \cM_\fg\\
G & \longmapsto G_\for.
\end{align*}
By definition, there is a factoring of  this map as
$$
\cM_p(n) \longr \cU(n) \longr \modfg
$$
through the open substack of formal groups of height at most $n$.
It is worth noting right away that the map $\cM_p(n) \to \cU(n)$
doesn't have a section. See Remark \ref{p-div-rem}.
\end{rem}

\begin{rem}[{\bf The geometry of $\cM_p(n)$}] This stack is something
of an mysterious object, despite years of work by many people.
Basic references include \cite{Messing}.
As an example of what is known, it has one geometric point
(i.e., isomorphism class of an algebraically closed field) for each integer
$h$, $1 \leq h \leq n$, given by the p-divisible group
$$
G_ h= \Ga_h \times (\ZZ/p^\infty)^{n-h}.
$$
Here  $\Ga_h$ is a formal group of height $h$ and 
$\ZZ/p^\infty$ is the colimit of the \'etale group schemes
$$
\ZZ/p^n = \Spec(\FF[x]/(x^{p^n}-x)).
$$
The morphism $\cM_p(n) \to \cU(n)$ is then surjective
on geometric points, but it is far from being an isomorphism.
For example, the automorphism group of $G_h$ is
$$
\Aut(\Ga_h) \times \mathrm{Gl}_{n-h}(\ZZ_p).
$$
\end{rem}

\begin{rem}\label{isitrep}
The morphism $\cM_p(n) \to \cM_\fg$ is {\it not} representable. This
follows from the statement about automorphisms
in the previous remark, but let's go into some detail. Consider
the two-category pull-back
$$
\xymatrix{
P_H \rto\dto & \cM_\fg(n)\dto\\
\Spec(R) \rto_H \rto & \modfg
}
$$
where $H$ is a formal group. By definition, $P_H$ is the functor
which assigns to each commutative ring $A$, the groupoid
of triples $(f,G,\phi)$
where $f:R \to A$ is a ring map, $G$ is a $p$-divisible group of
height $n$ over $A$ and $\phi:G_\for \to f^\ast H$ is an isomorphism.
Put another way, $P_H(A)$ is the  groupoid sheaf of extensions
\begin{equation}\label{sespb}
0 \to f^\ast H \to G \to G_\et \to 0
\end{equation}
over $A$. Isomorphisms fix $f^\ast H$, but not $G_\et$. If $P_H$
was actually equivalent to a scheme, then a short exact sequence of
the form \ref{sespb} would have no automorphisms, but this is
evidently not the case. To be specific, let $R=\bFF$ be the algebraic
closure of $\FF_p$ and let $\Ga_h$ be a formal group of height
$h$, $1 \leq h < n$ over $\bFF$. Then there is a split extension
$$
0 \to \Ga_h \to \Ga_h \times (\ZZ/p^\infty)^{n-h} \to (\ZZ/p^\infty)^{n-h} \to 0.
$$
There are no maps from $(\ZZ/p^\infty)^{n-h}$ to $\Ga_n$; therefore,
the automorphisms of this extension are $\mathrm{Gl}_{n-h}(\ZZ_p)$.
\end{rem}

We now can state Lurie's realization result. See \cite{Lurie}. 
Since we are working over $\ZZ_p$, one must take care with 
the hypotheses of here: the notions of algebraic stack
and \'etale must be the appropriate notions
over $\Spf(\ZZ_p)$.

\begin{thm}[{\bf Lurie}]\label{mainthm} Let $\cM$ be an algebraic stack
equipped with an \'etale morphism
$$
\cM \longr \cM_p(n).
$$
Then the realization problem for the composition
$$
\cM \longr \cM_p(n) \longr \cM_\fg
$$
has a canonical solution; that is, the space of all solutions 
is connected and has a preferred basepoint.
\end{thm}

\begin{rem}\label{john-wil}This theorem directly confronts the conundrum
of Example \ref{LEFT-realize}. Specifically, we can use this result to realize a
Landweber exact theory as an $E_\infty$ ring spectrum only if
the associated formal group is the formal part of $p$-divisible
group. This is a strong hypothesis. It works, for example, for $p$-complete
$K$-theory, but not evidently for the $p$-completed analog of
Johnson-Wilson theories
$E(n)_\ast$. See Remark \ref{p-div-rem}.4.

There is a deeper point, which I have put off discussing until
now. The homotopy groups $E_\ast$ of an $E_\infty$-ring spectum
$E$ support more structure
than that of a graded commutative ring. In particular, the operad
action maps
$$
(E\Sigma_n)_+ \wedge_{\Sigma_n} E^{\wedge n} \to E
$$
have induced maps in homotopy, which give rise to power operations
in $E_\ast$. This has a significant impact on the realization
problem --  for if $g:\Spec(R) \to \cM_\fg$ is a flat map classifying
a formal group $G$, there is no particular reason to suppose the
geometry of the formal group would specify the structure of the power
operations. However, if $g$ factors
$$
\xymatrix{
\Spec(R) \rto &\cM_p(n) \rto & \cM_\fg
}
$$
these power operations should be specified by the subgroup structure
of the $p$-divisible group.
\end{rem}

\begin{exam}[{\bf Serre-Tate theory}]\label{ST} As addendum to this theorem,
Lurie points out the morphism $\epsilon:\cM \to \cM_p(n)$ is \'etale if
it satisfies the Serre-Tate theorem; thus, for example, we recover
the Hopkins-Miller Theorem \ref{HM}, at least for
smooth elliptic curves.\footnote{In his monograph \cite{Lurie}, Lurie
suggests an argument for completing the proof of the Hopkins-Miller
Theorem, once we know the realization result for the open substack
$\mell$ of smooth elliptic curves.}

To state the Serre-Tate theorem we  need the language of deformation
theory. Let $\cM$ be a stack over $\cM_p(n)$ and $A_0/\FF$ be an $\cM$-object
over a field $\FF$, necessarily of characteristic $p$ since we are working
over $\Spf(\ZZ_p)$. Recall that an Artin local ring $(R,\mm)$ is a local ring
with nilpotent maximal ideal $\mm$. If $q:R \to \FF$
be a surjective morphism of rings, then a {\it deformation}
of $A_0$ to $R$ is an $\cM$-object $A$ over $R$ and a
pull-back diagram
$$
\xymatrix{
A_0 \rto \dto & A \dto\\
\Spec(\FF) \rto & \Spec(R).
}
$$
Deformations form a groupoid functor $\Def_\cM(\FF,A_0)$ on 
an appropriate category of Artin local rings.
The Serre-Tate theorem holds if the evident morphism
$$
\Def_\cM(\FF,A_0) \longr \Def_{\cM_p(n)}(\FF,\epsilon A_0)
$$
is an equivalence. This result holds for elliptic curves, but actually
in much wider generality. See \cite{Messing}.
\end{exam}

\begin{rem}[{\bf Deformations of $p$-divisible groups}]\label{Def} The deformation
theory of $p$-divisible groups and formal groups is well understood and
a simple application of Schlessinger's general theory \cite{Schlessinger}.
For formal groups, this is Lubin-Tate theory \cite{LT}. If $\Ga$ is
a formal group of height $n$ over a perfect field $\FF$, then Lubin-Tate
theory says that the groupoid-valued functor $\Def_{\cM_\fg}(\FF,\Ga)$
is discrete; that is, the  natural
map
$$
\Def_{\cM_\fg}(\FF,\Ga) \to \pi_0\Def_{\cM_\fg}(\FF,\Ga)
$$
is an equivalence. Furthermore, $\pi_0\Def_{\cM_\fg}(\FF,\Ga)$ is
pro-represented by a complete local ring $R(\FF,\Ga)$; that is,
there is a natural isomorphism
$$
\pi_0\Def_{\cM_\fg}(\FF,\Ga) \cong \Spf(R(\FF,\Ga)).
$$
A choice of $p$-typical coordinate for the universal deformation
of $\Ga$ over $R(\FF,\Ga)$ defines an isomorphism
$$
W(\FF)[[u_1,\ldots,u_{n-1}]] \cong R(\FF,\Ga)
$$
where $W(-)$ is the Witt vector functor.

A similar result holds for $p$-divisible groups. Let $G$ be a 
$p$-divisible group over an algebraically closed field $\FF$. Then we have 
split short exact sequence 
$$
0 \to G_\for \to G \to G_\et \to 0.
$$
Since $\FF$ is algebraically closed, there is an isomorphism
\begin{equation}\label{etalepart}
G_\et \cong (\ZZ/p^\infty)^{n-h}
\end{equation}
where $h$ is the height of $G_\for$.
Since $G_\et$ has a unique deformation up to isomorphism, by the definition
of \'etale, a choice of isomorphism \ref{etalepart} now identifies deformations
of $G$ as extensions
$$
0 \to H_\for \to H \to (\ZZ/p^\infty)^{n-h} \to 0
$$
where $H_\for$ is a deformation of $G_\for$. The exact sequence of sheaves
of groups
$$
0 \to \ZZ^{n-h} \to \QQ^{n-h} \to (\ZZ/p^\infty)^{n-h} \to 0
$$
now identifies the isomorphism class of extension as an element in
$$
\Hom(\ZZ^{n-h},H_\for).
$$
Thus we conclude that the groupoid-valued functor $\Def_{\cM_p(n)}(\FF,G)$
is discrete and $\pi_0\Def_{\cM_p(n)}(\FF,G)$ is pro-represented by
$$
R(\FF,\Ga)[[t_1,\cdots,t_{n-h}]] \cong W(\FF)[[u_1,\cdots,u_{h-1},t_1,\cdots,t_{n-h}]].
$$
Note that this is always a power series in $n-1$ variables. Similar results
hold for perfect fields by Galois descent.

Using this remark it is possible to give a local criterion for when a morphism
of stacks $\cM \to \cM_p(n)$ is \'etale. It is in this guise that Lurie's theorem
appears in \cite{BL}.
\end{rem}

\begin{rem}\label{the-proof}The proof of Theorem \ref{mainthm} has two large
steps. The first is to define and prove the existence
of an analog of the stack $\cM$ appropriate
for ``derived algebraic geometry'' -- which can be thought of as
algebraic geometry with $E_\infty$-ring spectra as the basic object.
This yields a stack $(\cM,\cO^\Top)$ where the structure sheaf is
now a sheaf of $E_\infty$-ring spectra. The second is to show that the
resulting algebraic object $(\cM,\cO^\Top)$ is the realization
required; that is, to construct an isomorphism $(\cM,\cO^\Top_\ast) \cong
(\cM,\cO_{\cM\ast})$. For this, there must be some homotopy theoretic input;
this is the local Hopkins-Miller theorem. This says that the Lubin-Tate
theory $E(\FF,\Ga)$ obtained from the deformations of a height $n$-formal
group over a perfect field $\FF$ is an $E_\infty$-ring spectrum and
the space of all $E_\infty$-structures is contractible. See \cite{GHbig}.
\end{rem}


\begin{thebibliography}{99}

\bibitem{Adams} J.\ F.\ Adams, \emph{Stable homotopy and generalised
cohomology}, University of Chicago Press, Chicago, 1974.

\bibitem{Bauer} Bauer, Tilman, ``Computation of the homotopy of thes
spectrum tmf'', Geometry and Topology Monographs 13 {\it Groups,
homotopy and configuration spaces (Tokyo 2005)}, 2008, 11--40.

\bibitem{BL} Behrens, M. and Lawson, T., `` Topological automorphic forms'',
available at {\tt http://front.math.ucdavis.edu/0702.5719}.

\bibitem{BTurner} Butowiez, Jean-Yves and Turner, Paul,
``Unstable multiplicative cohomology operations'', {\it Q. J. Math.},
51 (2000) no. 4, 437--449.

\bibitem{DelMod} Deligne, P., ``Courbes elliptiques: formulaire d'apr\`es
{J}. {T}ate'', {\it Modular functions of one variable, IV (Proc. Internat.
Summer School, Univ. Antwerp, Antwerp, 1972)}, 53--73, Lecture Notes
in Math., Vol. 476, Springer, Berlin 1975.

\bibitem{DM} Deligne, P. and Mumford, D., ``The irreducibility of the
space of curves of given genus'', {\it Inst. Hautes \'Etudes Sci. Publ. Math.},
36 (1969), 75--109.

\bibitem{DR} Deligne, P. and Rapoport, M., ``Les sch\'emas de
modules de courbes elliptiques'', {\it Modular functions of one
variable, {II} ({P}roc. {I}nternat.  {S}ummer {S}chool, {U}niv.
{A}ntwerp, {A}ntwerp, 1972)}, 143--316, Lecture Notes in Math., Vol. 349,
Springer, Berlin, 1973.

\bibitem{DG}  Demazure, Michel and Gabriel, Pierre,
{\it Groupes alg\'ebriques. {T}ome {I}: {G}\'eom\'etrie
alg\'ebrique, g\'en\'eralit\'es, groupes commutatifs},
Avec un appendice {\it Corps de classes local}\ par Michiel
Hazewinkel, Masson \& Cie, \'Editeur, Paris, 1970.

\bibitem{modfg} Goerss, P.G.., ``Quasi-coherent sheaves on the
moduli stack of formal groups'' {\tt http://front.math.ucdavis.edu/0802.0996}.

\bibitem{Newton} Goerss, P.G., ``(Pre-)sheaves of
ring spectra over the moduli stack of
formal group laws'' {\it Axiomatic, enriched and motivic
homotopy theory (NATO Advanced Study Institute)},  Kluwer Series C:
Mathematical and Physical Sciences 131, 101-131, Kluwer Academic Publishers,
Dordecht, 2004.

\bibitem{GHbig} Goerss, P.G. and Hopkins, M.J., ``Moduli spaces
of commutative ring
spectra'', {\it Structured ring spectra}, London Math. Soc. Lecture Note Ser.,
315, 151-200, Cambridge Univ. Press, Cambridge, 2004.
\smallskip

\bibitem{Hart} Hartshorne, Robin, {\it Algebraic geometry}, Graduate Texts
in Mathematics, No. 52, Springer-Verlag, New York, 1977.

\bibitem{HL} Hill, Michael and Lawson, Tyler, ``Automorphic forms and cohomology theories on Shimura curves of small discriminant'', {\bf arXiv:0902.2603}.

\bibitem{Hollander} Hollander, Sharon, ``Geometric criteria
for Landweber exactness'', preprint 2008.

\bibitem{ICM} Hopkins, M. J., ``Algebraic topology and modular forms'',
{\it Proceedings of the {I}nternational {C}ongress of
{M}athematicians, {V}ol. {I} ({B}eijing, 2002)}, 291-317,
Higher Ed. Press, Beijing, 2002.

\bibitem{HG}  M. J. Hopkins and B. H. Gross, ``The rigid analytic
period mapping, {L}ubin-{T}ate space, and stable homotopy theory'',
{\it Bull. Amer. Math. Soc. (N.S.)}, 30 (1994) No. 1, 76--86.

\bibitem{HopMah} Hopkins, M.J. and Mahowald, M.E, ``From elliptic
curves to homotopy theory'', manuscript available at
{\tt http://hopf.math.purdue.edu/}.

\bibitem{HS}  Hovey, Mark and Strickland, Neil P.,
``Morava {$K$}-theories and localisation'',
{\it Mem. Amer. Math. Soc.}, 139 (1999), no. 666.

\bibitem{Igusa} Igusa, Jun-ichi, ``Fibre systems of {J}acobian varieties.
{III}. {F}ibre systems of elliptic curves'', {\it Amer. J. Math.}, 81
(1959), 453--476.

\bibitem{Kashiwabara} Kashiwabara, Takuji, ``Hopf rings and unstable
operations'', {\it J. Pure Appl. Algebra}, 94 (1994), no.2, 183--193.

\bibitem{KM} Katz, Nicholas M. and Mazur, Barry,
{\it Arithmetic moduli of elliptic curves},
Annals of Mathematics Studies, 108,
Princeton University Press,
Princeton, NJ, 1985.

\bibitem{Kottwitz} Kottwitz, Robert E., ``Points on some {S}himura
varieties over finite fields'', {\it J. Amer. Math. Soc.}, 5 (1992) no. 2,
373-444.

\bibitem{Land} Landweber, P. S.,
``Homological properties of comodules over {$M{\rm U}\sb\ast
(M{\rm U})$}\ and {BP{$\sb\ast $}}({BP})'', {\it Amer. J. Math.},
98 (1976), no. 3, 591-610.

\bibitem{Laumon}  Laumon, G{\'e}rard and Moret-Bailly, Laurent,
{\it Champs alg\'ebriques},
{Ergebnisse der Mathematik und ihrer Grenzgebiete} 39, Springer-Verlag, Berlin, 2000.

\bibitem{Laz}  Lazard, Michel,
``Sur les groupes de {L}ie formels \`a un param\`etre'',
{\it Bull. Soc. Math. France}, 83 (1955), 251-274.

\bibitem{LT} J. Lubin and J. Tate, ``Formal moduli for
one-parameter formal {L}ie groups, {\it Bull. Soc. Math. France}, 94 (1966)
49--59.

\bibitem{Lurie} Lurie, J., ``Survey article on elliptic cohomology'',
manuscript 2007, available at {\tt http://www-math.mit.edu/\~{}lurie/}.

\bibitem{MahRez} Mahowald, Mark and Rezk, Charles,
``Brown-{C}omenetz duality and the {A}dams spectral sequence'',
{\it Amer. J. Math.}, 121 (1999), no. 6, 1153--1177.

\bibitem{Messing}  Messing, William,
{\it The crystals associated to {B}arsotti-{T}ate groups: with
applications to abelian schemes},
Lecture Notes in Mathematics, Vol. 264,
Springer-Verlag, Berlin, 1972.

\bibitem{LEFT} Miller, Haynes R., ``Sheaves, gradings, and the exact
functor theorem'' available from {\tt http://www-math.mit.edu/\~{ }hrm/papers/papers.html}

\bibitem{Neumann} Naumann, Niko, ``The stack of formal groups in stable
homotopy theory'', {\it Adv. Math.}, 215 (2007) no. 2, 569-600.

\bibitem{Rav} Ravenel, D.C., {\it Complex cobordism and stable
homotopy groups of spheres}, Academic Press Inc., Orlando, FL, 1986.

\bibitem{RavOrange} Ravenel, D.C., {\it Nilpotence and periodicity in
stable homotopy theory}, (Appendix C by Jeff Smith), Princeton University
Press, Princeton, NJ, 1992.

\bibitem{512} Rezk, Charles, ``Supplementary notes for Math 512'';
notes of a course on topological modular forms
given at Northwestern University in 2001, 
available at {\tt http://www.math.uiuc.edu/\~{}rezk/papers.html}.

\bibitem{Schlessinger} Schlessinger, Michael, ``Functors of {A}rtin rings'',
{\it Trans. Amer. Math. Soc.}, 130 (1968), 208-222.

\bibitem{Silver} Silverman, Joseph H., {\it The arithmetic of elliptic curves},
Graduate Texts in Mathematics 106, Springer-Verlag, New York 1986.

\bibitem{Smith} Smithling, Brian, ``On the moduli stack of commutative,
1-parameter formal Lie groups'', Thesis, University of Chicago, 2007.

\end{thebibliography}
\end{document}